\documentclass[12pt]{amsart}

\usepackage{amsmath,amssymb,amscd,amsthm,fullpage}
\usepackage[mathscr]{eucal}

\newtheorem{theorem}{Theorem}[section]
\newtheorem{proposition}[theorem]{Proposition}

\newtheorem{lemma}[theorem]{Lemma}

\theoremstyle{definition}

\newtheorem{definition}{Definition}[section]

\newcommand{\C}{\mathbb{C}}

\newcommand{\Z}{\mathbb{Z}}

\newcommand{\Q}{\mathbb{Q}}
\newcommand{\eps}{\varepsilon}
\newcommand{\p}{\partial}
\newcommand{\half}{\tfrac{1}{2}}
\newcommand{\CP}{\mathbb{CP}}
\renewcommand{\*}{\cdot}
\newcommand{\E}{\mathsf{E}}
\newcommand{\PP}{\mathsf{P}}
\newcommand{\<}{\langle}
\renewcommand{\>}{\rangle}

\DeclareMathOperator{\res}{res}

\renewcommand{\]}{{]\!]}}
\renewcommand{\[}{{[\![}}

\newcommand{\CO}{\mathcal{O}}
\renewcommand{\o}{\otimes}
\renewcommand{\)}{{)\!)}}
\renewcommand{\(}{{(\!(}}
\newcommand{\bull}{\bullet}
\renewcommand{\AA}{\mathcal{A}}

\newcommand{\CC}{\mathcal{C}}
\newcommand{\tAA}{\widetilde{\AA}}
\newcommand{\tII}{\widetilde{\II}}
\renewcommand{\d}{\delta}
\newcommand{\dbar}{\bar\d}
\renewcommand{\L}{\mathsf{L}}
\newcommand{\Lbar}{\bar{L}}
\newcommand{\Bbar}{\bar{B}}
\newcommand{\abar}{\bar{a}}
\newcommand{\pbar}{\bar{p}}
\newcommand{\qbar}{\bar{q}}
\newcommand{\vbar}{\bar{v}}
\newcommand{\xbar}{\bar{x}}
\newcommand{\z}{z}
\newcommand{\zbar}{\bar{\z}}
\renewcommand{\t}{\nu}
\newcommand{\II}{\mathcal{I}}
\newcommand{\LL}{\mathbb{L}}
\newcommand{\LLbar}{\bar{\LL}}
\newcommand{\CF}{\mathcal{F}}
\newcommand{\TT}{\mathbb{T}}
\newcommand{\s}[2]{\genfrac{[}{]}{0pt}{}{#1}{#2}}
\newcommand{\Mbar}{\overline{\mathcal{M}}}
\newcommand{\talpha}{\tilde{\alpha}}

\begin{document}

\title{The equivariant Toda lattice, I}

\author{Ezra Getzler}

\address{Department of Mathematics, Northwestern University, Evanston,
Illinois}

\maketitle

The Toda lattice is an infinite dimensional dynamical system of
commuting flows
$$
(\p,\d_n,\dbar_n \mid n>0) ,
$$
acting on functions $(q,a_k,\abar_k\mid k>0)$ defined on a one-dimensional
lattice. In the limit of small lattice spacing $\eps$, which is all that
will concern us here (Takasaki and Takebe \cite{TT}), the functions
$(q,a_k,\abar_k)$ become functions of a real parameter $x$, and the role of
translation by one unit of the lattice is taken by the operator
$\E=e^{\eps\p}$, where $\p$ is the infinitesimal generator of translations
in $x$.

The derivations $\d_1$ and $\dbar_1$ act on the variables $q$, $a_1$ and
$\abar_1$ by the formulas
\begin{align} \label{toda}
\d_1\abar_1 &= \dbar_1a_1 = \nabla q , &
\d_1\log q &= \nabla a_1 , & \dbar_1\log q &= \nabla\abar_1 ,
\end{align}
where $\nabla:\AA\to\AA$ is the infinite-order differential operator
\begin{align*}
\nabla &= \eps^{-1}\bigl(\E^{1/2}-\E^{-1/2}\bigr) =
\sum_{k=0}^\infty \frac{\eps^{2k}\p^{2k+1}}{2^{2k}(2k+1)!} \\
&= \p + \tfrac{1}{24} \, \eps^2 \, \p^3 + O(\eps^4) .
\end{align*}
These formulas imply the \textbf{Toda equation}:
\begin{equation} \label{Toda}
\d_1\dbar_1\log q = \nabla^2 q .
\end{equation}
The higher Toda flows are symmetries of this equation.

An abstract mathematical formulation of the Toda lattice is obtained by
realizing the derivations $(\p,\d_n,\dbar_n)$ on the free differential
algebra $\AA = \Q_\eps\{q,a_k,\abar_k\mid k>0\}$, defined over the ring
$\Q_\eps=\Q\[\eps\]$. Reductions of the Toda lattice are differential
ideals $\II$ in $\AA$ invariant under conjugation and closed under the
derivations $\d_n$ and $\dbar_n$. For example, the \textbf{Toda chain} is
described by the differential ideal with generators
$$
\{ a_1 - \abar_1 , a_2 - q , \abar_2 - q , a_k , \abar_k \mid k>2 \} .
$$

In this paper, we study a reduction of the Toda lattice, which we call the
equivariant Toda lattice. If $\t$ is a formal parameter, this reduction is
defined by the following constraint on the Lax operator:
\begin{equation} \label{w}
(\d_1 - \dbar_1)L = \t \p L .
\end{equation}
We prove that the corresponding differential ideal $\II_\t\subset\AA[\t]$
is isomorphic to the differential algebra
$$
\Q_{\eps,\t}\{q,v,\vbar\}/(\t\p q-\nabla(v-\vbar)) \o_{\Q_{\eps,\t}}
\Q_{\eps,\t}[\z_k,\zbar_k\mid k>0] ,
$$
where $\Q_{\eps,\t}=\Q_\eps[\t]$, $v$ and $\vbar$ are the images of
$a_1$ and $\abar_1\in\AA$, and $z_k$ and $\zbar_k$ are constants of
motion, which may be defined by the following equation:
$$
\biggl( L - \t + \sum_{k=1}^\infty z_k L^{-k} \biggr) \frac{\p L}{\p v} = L
.
$$

The equivariant Toda lattice may be used to describe the equivariant
Gromov-Witten invariants of $\CP^1$. Let $\TT$ be the multiplicative group
of $\C$, and let $X$ be a topological space with an action of $\TT$. The
equivariant cohomology $H^*_\TT(X,\Z)$ of $X$ is a module over the graded
ring $H^\bull_\TT(\ast,\Z)\cong\Z[\t]$, where $\t\in H^2_\TT(\ast,\Z)$. The
equivariant cohomology $H^\bull_\TT(\CP^1,\Z)$ of the projective line
$\CP^1$ admits a presentation
$$
H^\bull_\TT(\CP^1,\Z) \cong \Z[H,\t]/(H(H-\t)) ,
$$
where $H$ is the equivariant Chern class $c_1(\CO(1))\in H_\TT^2(\CP^1,\Z)$.

Denote the $k$th descendants of the cohomology classes $1$ and $H$ in
Gromov-Witten theory by $\tau_{k,P}$ and $\tau_{k,Q}$ respectively;
also, abbreviate $\tau_{0,P}$ and $\tau_{0,Q}$ to $P$ and $Q$. The
genus $0$ equivariant Gromov-Witten invariants of $\CP^1$ are
integrals over the moduli space $\Mbar_{0,n}(\CP^1)$ of stable maps of
genus $g$ with $n$ marked points:
$$
\<\tau_{k_1,P}\dots\tau_{k_m,P}\tau_{\ell_1,Q}\dots\tau_{\ell_n,Q}\>_g
\in H^\bull_\TT(\ast,\Q) \cong \Q[\t] .
$$
The large phase space is the formal affine space with
coordinates $\{s_k,t_k\mid k\ge0\}$. The genus~$g$ Gromov-Witten potential
$\CF_g$ of $\CP^1$ is the generating function on the large phase space given
by the formula
$$
\CF_g = \sum_{m,n=0}^\infty \frac{1}{m!\,n!}
\sum_{\substack{k_1,\dots,k_m \\ \ell_1,\dots,\ell_n}} s_{k_1}\dots s_{k_m}
t_{\ell_1} \dots t_{\ell_n}
\<\tau_{k_1,P}\dots\tau_{k_m,P}\tau_{\ell_1,Q}\dots\tau_{\ell_n,Q}\>_g
.
$$
We may combine the Gromov-Witten potentials into a single generating
function by interpreting $\eps$ as a genus expansion parameter, and
writing
$$
\CF = \sum_{g=0}^\infty \eps^{2g} \CF_g .
$$

Based on explicit calculations using the topological recursion relations in
genus $0$ and $1$, Pan\-dha\-ri\-pan\-de conjectured \cite{rahul} that the
following equation holds for the total Gromov-Witten potential:
\begin{equation} \label{rahul}
\p_0\bar{\p}_0\CF = \exp(\nabla^2\CF) .
\end{equation}
Here, $\nabla=\eps^{-1}(\E^{1/2}-\E^{-1/2})$, where $\E=e^{\eps\p}$, and
$\p=\p/\p s_0$. This equation was proved recently by Okounkov and
Pandharipande \cite{OP}.

On applying the operator $\nabla^2$ to both sides of \eqref{rahul} and
identifying the vector fields $\p_0$ and $\bar{\p}_0$ with the Toda flows
$\delta_1$ and $\dbar_1$, we obtain the Toda equation \eqref{Toda} for
$q=\exp(\nabla^2\CF)$. Observe that $\p_1-\bar{\p}_1=\t\p$; this equation
is formally identical to the constraint defining the equivariant Toda
lattice.

The equivariant Toda lattice and the equivariant Gromov-Witten theory of
$\CP^1$ each involve sequences $\{\delta_n,\dbar_n\}$ and $\{\p_n=\p/\p
t_n,\bar{\p}_n=\p/\p t_n-\t\p/\p s_n\}$ of commuting derivations, in the
first case on the algebra $\CC$, and in the second case on functions on the
large phase space. These sequences of vector fields may be compared by
means of a morphism
$$
\AA[\t]/\II_\t \longrightarrow \Q_{\eps,\t}\[s_k,t_k\mid k\ge0\]
$$
of differential algebras which sends the generators $q$, $v$ and $\vbar$ to
$\exp(\nabla^2\CF)$, $\nabla\p_0\CF$ and $\nabla\bar{\p}_0\CF$, and the
constants $z_k$ and $\zbar_k$ to $0$. In fact, the following relationship
between these flows holds:
\begin{subequations}
\begin{align}
\sum_{k=0}^\infty z^{k+1} \p_k &= \sum_{n=1}^\infty
\frac{z^n\delta_n}{(1+z\t)(2+z\t)\dots(n+z\t)} , \label{conjecture} \\
\sum_{k=0}^\infty z^{k+1} \bar{\p}_k &= \sum_{n=1}^\infty
\frac{z^n\dbar_n}{(1-z\t)(2-z\t)\dots(n-z\t)} . \label{conjecturebar}
\end{align}
\end{subequations}
We conjectured this in a preprint of this paper, based on a proof of the
result in genus 0 (see Section 5), together with calculations in genus 1
for small values of $n$; it has recently been proved by Okounkov and
Pandharipande \cite{OP}. Thus, the equivariant Toda lattice yields a
description of the equivariant Gromov-Witten invariants of $\CP^1$ in terms
of a Lax operator whose coefficients are obtained by an explicit recursion.

In a sequel to this paper, we relate the equivariant Toda lattice to the
dressing operator formalism. Let $\log(L)=W\log(\Lambda)W^{-1}$ be the
logarithm of the Lax operator $L$, related to the operator $\ell=\eps(\p
W)W^{-1}$ by the formula
$$
\log(L) = \log(\Lambda) - \ell .
$$
Borrowing ideas of Carlet, Dubrovin and Zhang \cite{CDZ}, we show that
the equivariant Toda lattice may be characterized by the expansion
$$
\Lambda + v + q\Lambda^{-1} = L + \t\ell - \sum_{k=1}^\infty \frac{z_k}{k}
L^{-k} .
$$
In particular, the equation $(\d_1-\dbar_1)W=\t\p W$ is equivalent to the
vanishing of the coefficients $z_k$. Under the same hypothesis, we also
show that the equivariant Toda lattice is Hamiltonian; this gives a more
direct relationship between the results of Okounkov and Pandharipande
\cite{OP} and the original Toda conjecture (Eguchi and Yang \cite{EY},
Eguchi, Hori and Yang \cite{EHY}, Pandharipande \cite{P}, Getzler
\cite{toda}).

\subsection*{Acknowledgements}

This paper has its origin in discussions with T. Eguchi and
R. Pan\-dha\-ri\-pan\-de at the workshop on Duality in Mirror Symmetry at
the Institute for Theoretical Physics of the University of California,
Santa Barbara. It was written during a membership of the Institute for
Advanced Study, Princeton. We are grateful to R. Dijkgraaf, B. Dubrovin and
Y. Zhang for useful discussions.

The author is partially supported by the NSF under grant DMS-0072508,
and, through the Institute for Advanced Study, under grant DMS-9729992.

{\appendix

\section{Another formulation of the equivariant Toda conjecture}

In this appendix, we invert the relations of \eqref{conjecture} and
\eqref{conjecturebar}. Define the (unsigned) Stirling numbers (of the first
kind) $\s{n}{k}$ by the generating function
$$
\sum_{k=0}^n \s{n}{k} \, \t^k = \prod_{j=0}^{n-1} (\t+j) .
$$
\begin{proposition}
\begin{subequations}
\begin{align}
\d_n &= n \sum_{k=1}^n \t^{k-1} \, \s{n}{k} \, \p_{n-k} ,
\label{Conjecture} \\
\dbar_n &= n \sum_{k=1}^n (-\t)^{k-1} \, \s{n}{k} \, \bar{\p}_{n-k} .
\label{Conjecturebar}
\end{align}
\end{subequations}
\end{proposition}
\begin{proof}
We prove the equivalence of \eqref{conjecture} with \eqref{Conjecture}; the
equivalence of \eqref{conjecturebar} and \eqref{Conjecturebar} is
similar. Equation \eqref{conjecture} may be restated as saying that
$$
\p_k = \sum_{n=1}^{k+1} (-\t)^{k-n+1}
h_{k-n+1}\bigl(1,\tfrac{1}{2},\dots,\tfrac{1}{n}\bigr)
\frac{\d_n}{n!} ,
$$
where $h_\ell$ is the complete symmetric polynomial of degree
$\ell$. We wish to prove that
$$
\frac{\d_n}{n!} = \sum_{\ell=0}^{n-1} \t^{n-\ell-1}
e_{n-\ell-1}\bigl(1,\tfrac{1}{2},\dots,\tfrac{1}{n-1}\bigr)
\p_\ell ,
$$
where $e_\ell$ is the elementary symmetric polynomial of degree
$\ell$. In other words, we wish to prove that
$$
\sum_{\ell=0}^{n-1} \t^{n-\ell-1} (-\t)^{\ell-m+1}
e_{n-\ell-1}\bigl(1,\tfrac{1}{2},\dots,\tfrac{1}{n-1}\bigr)
h_{\ell-m+1}\bigl(1,\tfrac{1}{2},\dots,\tfrac{1}{m}\bigr) =
\d_{n,m} .
$$
This is clearly true if $n\le m$; thus, we have only to prove that the
left-hand side vanishes when $n>m$. In this case, it equals $\t^{n-m}$
times the coefficient of $\t^{n-m}$ in the generating function
$$
\prod_{j=1}^{n-1} (1+j\t) \* \prod_{j=1}^m (1+j\t)^{-1} =
\prod_{j=m+1}^{n-1} (1+j\t) ,
$$
which is a polynomial of degree $n-m-1$; hence, the coefficient in
question vanishes.
\end{proof}
}

\section{Difference operators}

In this section, we recall the mathematical structure underlying the
Toda lattice; this material is adapted from the fundamental papers of
Ueno and Takasaki \cite{UT} and Kupershmidt \cite{K}.

All of the commutative algebras which we consider in this paper carry
an involution $p\mapsto\pbar$, and all ideals which we consider are
closed under this involution. By a \textbf{differential algebra}, we
mean a commutative algebra with derivation $\p$ such that
$$
\p\pbar = \overline{\p p} .
$$
A \textbf{differential ideal} is an ideal closed under the action of the
differential $\p$. If $S$ is a subset of a differential algebra $\AA$,
denote the differential ideal generated by $S\cup\bar{S}$ by $(S)$,
where $\bar{S}=\{\xbar\mid x\in S\}$ is the conjugate of $S$.

If $\AA$ is a differential algebra and $S$ is a set, the free
differential algebra $\AA\{S\}$ generated by $S$ is the polynomial
algebra
$$
\AA[\p^nx,\p^n\xbar\mid x\in S, n\ge0] ,
$$
with differential $\p(\p^nx)=\p^{n+1}x$.

An \textbf{evolutionary derivation} $\delta$ of a differential algebra $\AA$ is
a derivation such that $[\p,\delta] = 0$.  The evolutionary
derivations form a Lie subalgebra of the Lie algebra of derivations of
$\AA$, with involution
$$
\dbar p = \overline{\d\pbar} .
$$

Let $\AA$ be a differential algebra over $\Q_\eps$, and let $q\in\AA$
be a regular element (that is, having no zero-divisors) such that
$\qbar=q$. The localization $q^{-1}\AA$ of $\AA$ is a filtered
differential algebra, with differential $\p(q^{-1})=-q^{-2}\p q$.  Let
$\Phi_\pm(\AA,q)$ be the associative algebras of \textbf{difference
operators}
\begin{align*}
\Phi_+(\AA,q) &= \biggl\{ \sum_{k=-\infty}^\infty p_k \, \Lambda^k
\biggm| \text{$p_k\in q^{-k}\AA$, $p_k=0$ for $k\ll0$} \biggr\} , \\
\Phi_-(\AA,q) &= \biggl\{ \sum_{k=-\infty}^\infty p_k \, \Lambda^k
\biggm| \text{$p_k\in\AA$, $p_k=0$ for $k\gg0$} \biggr\} ,
\end{align*}
with product
$$
\sum_i a_i \Lambda^i \* \sum_j b_j \Lambda^j = \sum_k \biggl( \sum_{i+j=k}
\bigl(\E^{-j/2}a_i\bigr)\bigl(\E^{i/2}b_j\bigr) \biggr) \Lambda^k .
$$
Note that $\Phi_-(\AA,q)$ is in fact independent of $q$.

Let $A\mapsto A_\pm$ be the projections on $\Phi_\pm(\AA,q)$ defined
by the formulas
\begin{align*}
\biggl( \sum_{k=-\infty}^\infty p_k\Lambda^k \biggr)_+ &= \sum_{k=0}^\infty
p_k\Lambda^k , &
\biggl( \sum_{k=-\infty}^\infty p_k\Lambda^k \biggr)_- &=
\sum_{k=-\infty}^{-1} p_k\Lambda^k .
\end{align*}
We see that $A=A_-+A_+$. Define the residue
$\res:\Phi_pm(\AA,q)\to\AA$ by the formula
$$
\res\biggl( \sum_{k=-\infty}^\infty p_k\Lambda^k \biggr) = p_0 .
$$

For $k\in\Z$, let $[k]$ be the isomorphism of $\AA$
$$
[k] = \frac{\E^{k/2}-\E^{-k/2}}{\E^{1/2}-\E^{-1/2}} =
\sum_{j=1}^k \E^{(k+1)/2-j} = k + O(\eps^2) .
$$
Define $q^{[k]}$ by the recursion
$$
q^{[k+1]} = \E^kq \* \E^{-1/2}q^{[k]} ,
$$
with initial condition $q^{[0]}=1$. The involution
$$
A = \sum_{k=-\infty}^\infty p_k \Lambda^k \mapsto \bar{A} =
\sum_{k=1}^\infty \pbar_k \, q^{[k]} \Lambda^{-k} + \pbar_0 +
\sum_{k=1}^\infty \pbar_{-k} \, q^{-[k]} \Lambda^k ,
$$
defines an anti-isomorphism between the algebras $\Phi_+(\AA,q)$ and
$\Phi_-(\AA,q)$.

\section{The Toda lattice}

To formulate the Toda lattice, we introduce the differential algebra
$$
\AA = \Q_\eps\{q,a_k\mid k>0\}/(q-\qbar) .
$$
Since the generator $a_1$ plays a special role in the theory, we denote
it by $v$. It will be useful to define the symbol $a_0$ to equal $1$.

The Lax operator of the Toda lattice is the difference operator
$$
L = \Lambda + \sum_{k=1}^\infty a_k \Lambda^{-k+1} \in \Phi_-(\AA,q) ;
$$
its conjugate $\Lbar$ is given by the formula
$$
\Lbar = q\Lambda^{-1} + \sum_{k=1}^\infty \abar_k q^{-[k-1]} \Lambda^{k-1}
\in \Phi_+(\AA,q) .
$$
Introduce elements $p_k(n)\in\AA$, defined for all $n\ge0$ and $k\in\Z$:
$$
L^n = \sum_{k=-\infty}^n p_k(n) \Lambda^k .
$$

To define the evolutionary derivation $\d_n$ on the generators $a_k$
of $\AA$, introduce the difference operator $B_n=L^n_+$, and impose
the Lax equation $\d_nL=\eps^{-1}[B_n,L]$. This equation means that
$$
\eps^{-1}[B_n,L] = \sum_{k=1}^\infty \d_na_k \, \Lambda^{-k+1} .
$$
In order for this to be meaningful, it must be shown that
the coefficient of $\Lambda^k$ in $[B_n,L]$ vanishes for $k>0$. This
follows from the identity $[L^n,L]=0$: we have
$$
[B_n,L] = [B_n,L] - [L^n,L] = - [L^n_-,L] ,
$$
and it is clear that the coefficient of $\Lambda^k$ in $[L^n_-,L]$
vanishes if $k>0$. We also see that $\d_n a_k$ equals the coefficient of
$\Lambda^{-k+1}$ in
$$
- \eps^{-1} \sum_{j=0}^k [ p_{j-k}(n) \Lambda^{j-k}, a_j \Lambda^{-j+1} ] ,
$$
hence that
\begin{align} \label{a}
\d_n a_k &= \nabla p_{-k}(n) + \eps^{-1} \sum_{j=1}^{k-1} \Bigl(
\E^{(1-j)/2} p_{j-k}(n) \, \E^{(k-j)/2}a_j - \E^{(j-1)/2} p_{j-k}(n) \,
\E^{(j-k)/2}a_j \Bigr) \\
&= \nabla p_{-k}(n) + \sum_{j=1}^{k-1} \Bigl( \E^{(1-j)/2} p_{j-k}(n) \,
\nabla [k-j] a_j - \E^{(j-k)/2}a_j \nabla [j-1] p_{j-k}(n) \Bigr) .
\notag
\end{align}
In particular, $\d_nv = \nabla p_{-1}(n)$.

To define $\d_n$ on the remaining generators $q$ and $\abar_k$ of
$\AA$, we impose the Lax equation $\d_n\Lbar=\eps^{-1}[B_n,\Lbar]$. In
particular, we see that
\begin{equation} \label{q}
\d_nq = q\nabla p_0(n) ,
\end{equation}
and hence that $\d_n q^{[k]} = q^{[k]} \nabla [k] p_0(n)$.
It also follows that
$$
\d_n(q^{-[k-1]}\abar_k) = q^{-[k-1]} \bigl(
\d_n\abar_k-\abar_k\nabla[k-1]p_0(n) \bigr)
$$
equals the coefficient of $\Lambda^{k-1}$ in
$$
\eps^{-1} \sum_{j=0}^k [ p_{k-j}(n) \Lambda^{k-j}, q^{-[j-1]} \abar_j
\Lambda^{j-1} ] ,
$$
hence that
\begin{align} \label{abar}
\d_n \abar_k &= \abar_k \nabla[k-1]p_0(n) \\
& \qquad + \eps^{-1} q^{[k-1]} \sum_{j=0}^k \Bigl( \E^{(1-j)/2} p_{k-j}(n) \,
\E^{(k-j)/2}q^{-[j-1]}\abar_j - \E^{(j-1)/2} p_{k-j}(n) \,
\E^{(j-k)/2}q^{-[j-1]}\abar_j \Bigr) \notag \\
&= \abar_k \nabla[k-1]p_0(n) \notag \\
& \qquad + \eps^{-1} \sum_{j=0}^k \Bigl( \E^{(1-j)/2}
\bigl( q^{[k-j]}p_{k-j}(n) \bigr) \, \E^{(k-j)/2}\abar_j
- \E^{(j-1)/2} \bigl( q^{[k-j]}p_{k-j}(n) \bigr) \, \E^{(j-k)/2}\abar_j
\Bigr) \notag \\
&= \nabla \bigl( q^{[k]}p_k(n) \bigr) \notag \\
&\qquad + \sum_{j=1}^{k-1} \Bigl( \E^{(1-j)/2} \bigl( q^{[k-j]}p_{k-j}(n)
\bigr) \, \nabla [k-j] \abar_j - \E^{(j-k)/2}\abar_j \nabla [j-1] \bigl(
q^{[k-j]}p_{k-j}(n) \bigr) \Bigr) . \notag
\end{align}
In particular, $\d_n\vbar = \nabla\bigl(qp_1(n)\bigr)$.

Denote by $\alpha:\AA\to\Q_\eps$ the homomorphism which sends the
generators $\{q,a_k,\abar_k\}$ of $\AA$ to $0$. By formulas \eqref{a},
\eqref{q} and \eqref{abar}, we see that $\d_nq$, $\d_na_k$ and
$\d_n\abar_k$ all lie in the ideal $(\p\AA)$ of $\AA$, and hence
\begin{equation} \label{deltaalpha}
\alpha\*\d_n=0 .
\end{equation}

We now recall the proof that the derivations $\d_m$ and $\d_n$
commute. The proof relies on the Zakharov-Shabat equation
\begin{equation} \label{ZS1}
\d_mB_n - \d_nB_m = \eps^{-1} [B_m,B_n] .
\end{equation}
To prove this equation, observe that
\begin{align*}
\d_mB_n &= (\d_mL^n)_+ = \eps^{-1} [B_m,L^n]_+ =
\eps^{-1} [B_m,B_n+L^n_-]_+ \\
&= \eps^{-1} [B_m,B_n] + \eps^{-1} [B_m,L^n_-]_+ .
\end{align*}
Since $[L^m,L^n]=0$, we also see that
$$
[B_m,L^n]_+ = [L^m-L^m_-,L^n]_+ = - [L^m_-,B_n]_+ .
$$
It follows that
$$
\d_mB_n - \d_nB_m = \eps^{-1} \bigl( [B_m,B_n] + [B_m,L^n_-]_+ \bigr) +
\eps^{-1}[L^n_-,B_m]_+ = \eps^{-1} [B_m,B_n] .
$$
From \eqref{ZS1}, we easily see that the derivations $\d_m$ and $\d_n$
commute:
\begin{align*}
[\d_m,\d_n]L &= \eps^{-1} \d_m[B_n,L] - \eps^{-1} \d_n[B_m,L] \\
&= \eps^{-1} [\d_mB_n,L] + \eps^{-2} [B_n,[B_m,L]] - \eps^{-1}
[\d_nB_m,L] - \eps^{-2} [B_m,[B_n,L]] \\
&= \eps^{-1} [\d_mB_n-\d_nB_m,L] + \eps^{-2} [B_n,[B_m,L]] - \eps^{-2}
[B_m,[B_n,L]] = 0 .
\end{align*}

The derivation $\dbar_n$ is defined to be the conjugate of $\d_n$, acting
on the generators of $\AA$ by the formulas
\begin{align*}
\dbar_nq &= \overline{\d_nq} , &
\dbar_na_k &= \overline{\d_n\abar_k} , &
\dbar_n\abar_k &= \overline{\d_na_k} .
\end{align*}
The following proposition establishes the Lax equation for this
derivation.
\begin{proposition}
Let $C_n=-\Lbar^n_-$; then $\dbar_nL=\eps^{-1}[C_n,L]$ and
$\dbar_n\Lbar=\eps^{-1}[C_n,\Lbar]$.
\end{proposition}
\begin{proof}
We have
\begin{align*}
\overline{\dbar_n\Lbar} &= \overline{\dbar_nq} \Lambda^{-1} +
\sum_{k=1}^\infty \overline{q^{-[k-1]} \bigl( \dbar_n\abar_k - \abar_k
q^{-[k-1]} \dbar_nq^{[k-1]} \bigr) \Lambda^{k-1}} \\
&= \sum_{k=1}^\infty \overline{\dbar_n\abar_k} \Lambda^{-k+1} + \nabla
p_0(n) \Lambda - \sum_{k=1}^\infty a_k \nabla[k-1]p_0(n) \Lambda^{-k+1} \\
&= \sum_{k=1}^\infty \d_na_k \Lambda^{-k+1} - [p_0(n),L]
= \eps^{-1} [B_n,L] - [p_0(n),L] \\
&= \eps^{-1} \overline{[C_n,\Lbar]} .
\end{align*}
A similar proof shows that
$\overline{\dbar_nL}=\eps^{-1}\overline{[C_n,L]}$.
\end{proof}

It is automatic that the derivations $\dbar_m$ and $\dbar_n$ commute,
since their conjugates do. To see that $\d_m$ commutes with $\dbar_n$,
we use the Zakharov-Shabat equation
\begin{equation} \label{ZS2}
\d_m\Bbar_n - \dbar_nB_m = \eps^{-1} [B_m,\Bbar_n] .
\end{equation}
This is proved by combining the equations
$$
\d_m\Bbar_n = (\d_m\Lbar^n)_- = \eps^{-1} [B_m,\Lbar^n]_- =
\eps^{-1} [B_m,\Bbar_n]_- ,
$$
and
$$
\dbar_mB_n = (\dbar_nL^m)_+ = \eps^{-1} [\Bbar_n,L^m]_+ =
\eps^{-1} [\Bbar_n,B_m]_+ .
$$
It follows from \eqref{ZS2} that $\d_m$ and $\dbar_n$ commute:
\begin{align*}
[\d_m,\dbar_n]L &= \eps^{-1} \d_m[\Bbar_n,L] - \eps^{-1} \dbar_n[B_m,L] \\
&= \eps^{-1} [\d_m\Bbar_n,L] + \eps^{-2} [\Bbar_n,[B_m,L]] - \eps^{-1}
[\dbar_nB_m,L] - \eps^{-2} [B_m,[\Bbar_n,L]] \\
&= \eps^{-1} [\d_m\Bbar_n-\dbar_nB_m,L] + \eps^{-2} [\Bbar_n,[B_m,L]] -
\eps^{-2} [B_m,[\Bbar_n,L]] = 0 . 
\end{align*}

\section{The equivariant Toda lattice}

Informally, a reduction of the Toda lattice is an invariant submanifold of
the configuration space fixed by the involution; we formalize this as
follows.
\begin{definition} \label{reduction}
A \textbf{reduction} of the Toda lattice is a differential ideal
$\II\subset\AA$ invariant under conjugation and preserved by the action of
the derivations $\d_n$ and $\dbar_n$.
\end{definition}

The simplest example of such a reduction is the \textbf{Toda chain},
defined by the constraint $L=\Lbar$; the associated differential ideal
$$
\II = ( v - \vbar , a_2-q , a_k \mid k>2 )
$$
is generated by the coefficients of $L-\Lbar$. To see that $\II$ is closed
under the action of the derivations $\d_n$, it suffices to observe that the
operator $L-\Lbar$ satisfies the Lax equation
$\d_n(L-\Lbar)=\eps^{-1}[B_n,L-\Lbar]$, and that the coefficients of
$\eps^{-1}[B_n,L-\Lbar]$ are contained in the differential ideal generated
by the coefficients of $L-\Lbar$.

The constraint $L=\Lbar$ is equivalent to the relation $\d_n=\dbar_n$ among
the Toda flows, for all $n$; in particular, the Toda equation \eqref{Toda}
becomes in this limit the equation $\d_1^2\log q=\nabla^2 q$.

In this paper, we study a reduction of the Toda lattice which is a
deformation of the Toda chain. Let $\AA[\t]$ be the extension of the
differential algebra $\AA$ by a variable $\t$, such that $\p\t=0$ and
$\bar{\t}=-\t$, and consider families of reductions of the Toda lattice
parametrized by $\AA[\t]$; that is, we consider differential ideals in
$\AA[\t]$ satisfying the conditions of Definition \ref{reduction}.

\begin{definition}
The \textbf{equivariant Toda lattice} is the reduction of the Toda lattice
defined over $\AA[\t]$ by the constraints
\begin{align} \label{equivariant}
(\d_1 - \dbar_1)L &= \t\p L , & (\d_1 - \dbar_1)\Lbar &= \t\p\Lbar .
\end{align}
Let $\II_\t$ be the differential ideal determining this reduction.
\end{definition}

Let $K$ be the difference operator
\begin{equation} \label{K}
K = B_1-C_1 = \Lambda+v+q\Lambda^{-1} .
\end{equation}
Substituting the Lax equations into the constraints \eqref{equivariant}, we
obtain an equivalent formulation of the equivariant Toda lattice: it is
characterized by the pair of equations
\begin{align} \label{Equivariant}
\eps^{-1} [K,L] &= \t\p L , & \eps^{-1}[K,\Lbar] &= \t\p\Lbar .
\end{align}
In other words, the differential ideal $\II_\t$ defining the equivariant
Toda lattice is generated by the coefficients of these equations.

The following theorem collects the main properties of the differential
ideal $\II_\t$.
\begin{theorem} \label{Main}
Let $\PP:\AA\to\AA$ be the infinite-order differential operator
\begin{align*}
\PP &= \frac{\p}{\nabla} = \sum_{g=0}^\infty
\frac{\eps^{2g}(2^{1-2g}-1)B_{2g}}{(2g)!}  \, \p^{2g} \\
&= 1 - \tfrac{1}{24} \, \eps^2 \, \p^2 + O(\eps^4) .
\end{align*}
The differential ideal $\II_\t$ defining the equivariant Toda lattice
equals $( y,\p \z_k\mid k>0)$, where
\begin{align*}
y &= q\nabla(v-\vbar) - \t\p q , \\
\z_k &= p_{-1}(k) - qp_1(k) - \t\PP p_0(k) .
\end{align*}
The differential algebra $\AA[\t]/\II_\t$ is isomorphic to
$$
\tAA = \Q_{\eps,\t}\{q,v,\z_k\mid k>0\}/(q-\qbar,y,\p\z_k\mid k>0) .
$$
The Toda flows $\d_n$ and $\dbar_n$ map $\z_k$ to $\II_\t$. 
\end{theorem}
\begin{proof}
Let $\tII_\t$ be the differential ideal $(y,\p\z_k\mid k>0)$. Define
elements $f_k(n)$, $g_k(n)\in\AA$ by the formulas
\begin{align*}
\eps^{-1}[K,L^n] - \t\p L^n &= \sum_{k=-\infty}^\infty f_k(n) \Lambda^{-k} , &
\eps^{-1}[K,\Lbar^n] - \t\p\Lbar^n &= \sum_{k=-\infty}^\infty g_k(n)
\Lambda^k .
\end{align*}
The differential ideal $\II_\t$ is generated by the coefficients
$f_k=f_k(1)$ and $g_k=g_k(1)$. The formulas
\begin{align*}
\eps^{-1}[K,L^n] - \t\p L^n &= \sum_{i=1}^n L^{i-1} \bigl( \eps^{-1}[K,L] -
\t\p L \bigr) L^{n-i} , \\
\eps^{-1}[K,\Lbar^n] - \t\p\Lbar^n &= \sum_{i=1}^n \Lbar^{i-1} \bigl(
\eps^{-1}[K,\Lbar] - \t\p\Lbar \bigr) \Lbar^{n-i} ,
\end{align*}
show that the coefficients $f_k(n)$ and $g_k(n)$ lie in $\II_\t$; hence
$\p\z_n=\PP f_0(n)$ and $\p\zbar_n=\PP g_0(n)$ do as well, showing that
$\tII_\t\subset\II_\t$. We wish to prove the equality of these two
differential ideals.

To do this, we show, by induction on $k$, that the coefficients $f_k$ lie
in $\tII_\t$. We have $f_k=0$ for $k<0$. If $f_j\in\tII_\t$ for $j<k$, we
see that
\begin{align*}
\nabla \z_{k+1} &= \res\bigl( \eps^{-1}[K,L^{k+1}] - \t\p L^{k+1} \bigr) =
\sum_{i=0}^k \res \bigl( L^i ( \eps^{-1}[K,L] - \t\p L ) L^{k-i} \bigr) \\
&\equiv [k+1] f_k \pmod{\tII_\t} ,
\end{align*}
hence $f_k\in\tII_\t$. A similar induction shows that $g_k\in\tII_\t$; this
induction starts with the fact that $g_{-1}=y$ lies in $\tII_\t$.

Since $\z_k-[k]a_{k+1}\in(q,v,a_2,\dots,a_k)$, we see that the
differential algebra $\AA/\II_\t$ is isomorphic to
$$
\tAA = \Q_{\eps,\t}\{q,v,\z_k\mid k>0\}/(q-\qbar,y,\p\z_k\mid k>0) .
$$

It remains to prove that $\d_n\z_k$ and $\dbar_n\z_k$ lie in $\II_\t$. By the
Zakharov-Shabat equations \eqref{ZS1} and \eqref{ZS2}, we see that
\begin{align*}
\d_nK &= \d_n(B_1-C_1) = (\d_1-\dbar_1)B_n + \eps^{-1}[B_n,B_1-C_1] \\
&= \t \p B_n + \eps^{-1} [B_n,K] .
\end{align*}
It follows that
\begin{align*}
\nabla\d_n \z_k &= \d_nf_0(k) = \d_n \res(\eps^{-1}[K,L^k]-\t\p L^k) \\
&= \res(\eps^{-1}[\d_nK,L^k]+\eps^{-1}[K,\d_nL^k]-\t\p\d_nL^k) \\
&= \eps^{-1} \res([\t\p B_n+\eps^{-1}[B_n,K],L^k] +
\eps^{-1}[K,[B_n,L^k]] - \t\p[B_n,L^k]) \\
&= \eps^{-1} \res([B_n,\eps^{-1}[K,L^k]-\t\p L^k]) \\
&= \nabla \sum_{j=1}^n [j] \bigl(p_j(n)f_j(k)\bigr) .
\end{align*}

The extension of $\alpha$ to a homomorphism from $\AA[\t]$ to
$\Q_{\eps,\t}$ continues to satisfy \eqref{deltaalpha}. It follows that
$\alpha(\d_n\z_k)=0$, hence we obtain an explicit equation for $\d_n\z_k$:
$$
\d_n\z_k = \sum_{j=1}^n [j] \bigl(p_j(n)f_j(k)\bigr) \in\II_\t .
$$
The proof that $\dbar_n\z_k\in\II_\t$ follows along the same lines.
\end{proof}

Let us calculate explicitly the coefficients $a_2$ and $a_3$ of the Lax
operator $L$ as elements of $\tAA$. Applying $\res:\Phi_-(\AA,q)\to\AA$ to
the equation $\eps^{-1}[K,L]=\t\p L$, we see that
$$
a_2 = q + \t \PP v + \z_1 .
$$
Taking the coefficient of $\Lambda^{-1}$ in the equation
$\eps^{-1}[K,L]=\t\p L$, we see that
$$
\nabla a_3 + (a_2-q) \nabla v = \t \p a_2 .
$$
\begin{lemma} \label{dP}
$\nabla f \, \PP g = \half \nabla \bigl( f [2]\PP g \bigr) - \half [2]
\bigl( f \p g \bigr)$
\end{lemma}
\begin{proof}
We have
$$
\nabla f \, \PP g = \eps^{-1} \E^{1/2} \Bigl( f \, \E^{-1/2}\PP g \bigr)
- \eps^{-1} \E^{-1/2} \bigl( f \, \E^{1/2}\PP g \bigr) .
$$
The result follows, since $\E^{\pm1/2}\PP=\half[2]\PP\pm\half\eps\p$.
\end{proof}

By this lemma,
$$
(a_2-q) \nabla v = \t \nabla v \, \PP v + \z_1 \nabla v = \t
\nabla \bigl( \half v [2]\PP v - \tfrac14 \PP [2]v^2 \bigr) + \z_1
\nabla v .
$$
It follows that
$$
a_3 = \t \bigl( \PP \bigl( \tfrac{1}{4} [2]v^2 + q \bigr) - \half v
[2]\PP v \bigr) + \t^2 \PP v - \z_1 v + \half \z_2 .
$$

This method of calculating the coefficients $a_k$ becomes cumbersome for
larger values of $k$: instead, it is better to use the recursion
\begin{equation} \label{recursez}
p_{-1}(n) = qp_1(n) + \t\PP p_0(n) + \z_n
\end{equation}
which is a consequence of Theorem \ref{Main}.

Let $\Psi$ be the algebra of difference operators
$$
\Psi = \{ A\in\Phi_-(\tAA,q) \mid \eps^{-1}[K,A] = \t\p A \} .
$$
Let $\LL\in\Psi$ be the Lax operator defined by the recursion
\begin{equation} \label{recurse}
p_{-1}(n) = qp_1(n) + \t\PP p_0(n) .
\end{equation}
This Lax operator plays a special role in the theory: the following lemma
shows that the algebra $\Psi$ may be identified with the commutative
algebra $\tAA_0\(\LL\)$, where
$$
\tAA_0 = \Q_{\eps,\t}[\z_k,\zbar_k\mid k>0]
$$
is the kernel of the derivation $\p:\tAA\to\tAA$.
\begin{lemma}
The homomorphism $\talpha:\tAA\to\tAA_0$ which sends the generators
$\p^nq$, $\p^nv$ and $\p^n\vbar$ of $\tAA$ to $0$ induces an isomorphism
between $\Psi$ and $\tAA_0\(\Lambda^{-1}\)$.
\end{lemma}
\begin{proof}
Since $\talpha(\LL)=\Lambda$, the map
$\talpha:\Psi\to\tAA_0\(\Lambda^{-1}\)$ is surjective. Suppose that
$A\in\Psi$ lies in the kernel of $\talpha$, and let $k$ be the smallest
integer such that the coefficient $x\in\tAA$ of $\Lambda^{-k}$ in $A$ is
nonzero. We have
$$
\eps^{-1}[K,A] - \t\p A = \nabla x \, \Lambda^{1-k} + O(\Lambda^{-k}) ,
$$
hence $x\in\tAA_0$. In this way, we see that
$\talpha:\Psi\to\tAA_0\(\Lambda^{-1}\)$ is injective.
\end{proof}

\begin{theorem}
The evolutionary derivation $e=\p_v+\p_{\vbar}$ of $\tAA$ preserves $\Psi$,
and $e(L)$ satisfies the formula
$$
\biggl( L - \t + \sum_{k=1}^\infty z_k L^{-k} \biggr) e(L) = L .
$$
\end{theorem}
\begin{proof}
If $A\in\Psi$, we have
$$
\eps^{-1}[K,e(A)]-\t\p e(A) = e(\eps^{-1}[K,A]-\t\p A) - [e(K),A] = 0 ,
$$
since $e(K)=1$. This shows that $e$ preserves $\Psi$.

For $n>0$, we have by \eqref{recurse} that $\talpha(e(p_{-1}(n))) = \t
\talpha(e(p_0(n)))$, or equivalently,
$$
\oint (\LL-\t) e(\LL^n) \frac{d\LL}{\LL} = 0 .
$$
Since $\Psi$ is a commutative algebra, $e(\LL)$ commutes with $\LL$, hence
$e(\LL^n)=n\LL^{n-1}e(\LL)$, and
$$
\oint (\LL-\t) \LL^{n-2} e(\LL) \, d\LL = 0 , \quad n>0 .
$$
This shows that the coefficient of $\LL^{-k}$ in $(\LL-\t)e(\LL)$ vanishes,
hence $(\LL-\t)e(\LL)=\LL$.

Since $L/e(L)$ lies in $\Psi$, there is an expansion
$$
\frac{L}{e(L)} = L - \t + \frac{1}{2\pi i} \sum_{n=0}^\infty L^{-n-1} \oint
L^n \frac{dL}{e(L)} .
$$
(The constant term is determined by the fact that
$e(L)=1+\t\Lambda^{-1}+O(\Lambda^{-2})$.) We have
$$
\frac{dL}{e(L)} = \frac{d\LL}{e(\LL)} = (\LL-\t) \frac{d\LL}{\LL} ,
$$
hence
$$
\frac{1}{2\pi i} \oint L^n \frac{dL}{e(L)} = \frac{1}{2\pi i} \oint
L^n \, (\LL-\t) \frac{d\LL}{\LL} .
$$
It follows from the recursion \eqref{recursez} that
$$
z_n = \talpha(p_{-1}(n)) - \t\talpha(p_0(n)) = \frac{1}{2\pi i} \oint
L^n \, (\LL-\t) \frac{d\LL}{\LL} ,
$$
and the theorem follows.
\end{proof}

\section{The dispersionless limit of the equivariant Toda lattice}
\label{genus0}

In this section, we consider the dispersionless limit of the
equivariant Toda lattice, in which $\eps\to0$; we only consider the
case in which the constants of motion $\z_k$ are setto $0$. If
$A\in\Phi_\pm(\tAA,q)$, we write
$$
A_0 = \lim_{\eps\to0} A .
$$
In the dispersionless limit, the algebra $\Phi_-(\tAA,q)$ degenerates to
the commutative algebra $\tAA\(\Lambda^{-1}\)$, and the leading order in
the commutator is the Poisson bracket
$$
\{A_0,B_0\} = \lim_{\eps\to0} \eps^{-1} [ A , B ] = (\Lambda\p_\Lambda
A_0)\p B_0 - \p A_0 (\Lambda\p_\Lambda B_0) .
$$
It is not hard to write down explicit formulas for the Lax operator $\LL$
of the equivariant Toda lattice and its conjugate $\LLbar$ in the
dispersionless limit.
\begin{theorem} \label{LL0}
We have
\begin{align*}
\LL_0 &= K_0 + \t \sum_{n=0}^\infty \Bigl(\frac{\t}{\Lambda}\Bigr)^n
\sum_{k=0}^n (-1)^{n-k} \s{n}{k} \frac{\log(K_0/\Lambda)^{n-k+1}}{(n-k+1)!}
\Bigl(\frac{K_0}{\Lambda}\Bigr)^{-n} \in \Q[\t,q,v]\(\Lambda^{-1}\)
\intertext{and}
\LLbar_0 &= K_0 - \t \sum_{n=0}^\infty \Bigl(\frac{-\t\Lambda}{q}\Bigr)^n
\sum_{k=0}^n (-1)^{n-k} \s{n}{k} \frac{\log(\Lambda K_0)^{n-k+1}}{(n-k+1)!}
\Bigl( \frac{\Lambda K_0}{q} \Bigr)^{-n} \in \Q[\t,q,\vbar]\(\Lambda^{-1}\) .
\end{align*}
\end{theorem}
\begin{proof}
Denote by $L_0$ the expression which we wish to prove equals
$\LL_0$. It is clear that $\beta(L)=\Lambda$, hence it suffices to
prove the equation
$$
\{K_0,L_0\} = \t \p K_0 ,
$$
which is the dispersionless limit of the equation $\eps^{-1}[K,\LL]=\t\p\LL$.

Since $\{K_0,\log(K_0/\Lambda)\}=\p K_0$, we have
$$
\{K_0,L_0\} = \t \p K_0 \sum_{n=0}^\infty \Bigl(\frac{\t}{\Lambda}\Bigr)^n
\sum_{k=0}^n (-1)^{n-k} \s{n}{k} \frac{\log(K_0/\Lambda)^{n-k}}{(n-k)!}
\Bigl(\frac{K_0}{\Lambda}\Bigr)^{-n} .
$$
Since $\p\log(K_0/\Lambda)=K_0^{-1}\p K_0$, we have
$$
\p L_0 = \p K_0 \biggl( 1 + \sum_{n=0}^\infty
\Bigl(\frac{\t}{\Lambda}\Bigr)^n
\sum_{k=0}^n (-1)^{n-k} \s{n}{k} \biggl(
\frac{\log(K_0/\Lambda)^{n-k}}{(n-k)!} -
\frac{n\log(K_0/\Lambda)^{n-k+1}}{(n-k+1)!} \biggr)
\Bigl( \frac{K_0}{\Lambda} \Bigr)^{-n-1} \biggr) .
$$
The equation $\{K_0,L_0\}=\t\p L_0$ follows from the recursion
$\s{n}{k} = (n-1) \s{n-1}{k} + \s{n-1}{k-1}$.

Define $e^{\circ j}(L_0)$ by induction: $e^{\circ 0}(L_0)=L_0$ and
$e^{\circ(j+1)}(L_0)=e(e^{\circ j}(L_0))$. Then
$$
e^{\circ j}(L_0) = \delta_{j,0} \, K_0 + \t^{1-j} \sum_{n=0}^\infty
\Bigl(\frac{\t}{\Lambda}\Bigr)^n \sum_{k=0}^{n-j+1} (-1)^{n-k-j} \s{n}{k}
\frac{\log(K_0/\Lambda)^{n-k-j+1}}{(n-k-j+1)!}
\Bigl(\frac{K_0}{\Lambda}\Bigr)^{-n} .
$$
This formula is proved by induction on $j$, using the
formulas $e(K_0)=1$ and $e(\log(K_0/\Lambda))=K_0^{-1}$.

There is an embedding of the differential algebra $\tAA$ in the differential
algebra
$$
\tAA\{u\}/(\p q-q\p u) \cong \Q_{\eps,\t}[q]\{u,v\} ,
$$
given by mapping $\vbar$ to $v-\t\PP u$. In the dispersionless limit,
this embedding maps $\vbar$ to $v-\t u$. We will prove the formula
for $\LLbar_0$ by working with Laurent series in this larger algebra.

The Laurent series $\LLbar_0$ is obtained from $\LL_0$ by replacing $v$ by
$v-\t u$, $\Lambda$ by $q/\Lambda$, and $\t$ by $-\t$. Let
$\widetilde{\LL}{}_0$ be the result of substituting $v-\t u$ for $v$ in
$\LL_0$; it is given by the formula
\begin{align*}
\widetilde{\LL}{}_0 &= \sum_{j=0}^\infty \frac{(-\t u)^j}{j!} e^{\circ j}(L) \\
&= K_0 + \t \sum_{n=0}^\infty \Bigl(\frac{\t}{\Lambda}\Bigr)^n \sum_{j=0}^{n+1}
\sum_{k=0}^{n-j+1} (-1)^{n-k} \s{n}{k}
\frac{u^j\log(K_0/\Lambda)^{n-k-j+1}}{j!(n-k-j+1)!}
\Bigl(\frac{K_0}{\Lambda}\Bigr)^{-n} \\
&= K_0 + \t \sum_{n=0}^\infty \Bigl(\frac{\t}{\Lambda}\Bigr)^n
\sum_{k=0}^{n+1} (-1)^{n-k} \s{n}{k}
\frac{\log(qK_0/\Lambda)^{n-k+1}}{(n-k+1)!}
\Bigl(\frac{K_0}{\Lambda}\Bigr)^{-n} .
\end{align*}
We obtain $\LLbar_0$ on substituting $-\t$ for $\t$ and $\Lambda$ for
$q/\Lambda$.
\end{proof}

Define generating functions $\pi_k(z)\in\Q[\t,q,v]\(z\)$ by the formula
$$
\LL_0(z) = \sum_{n=-\infty}^\infty \frac{z^n}{[n]!} \LL_0^n =
\sum_{k=-\infty}^\infty \pi_k(z) \Lambda^k ,
$$
where $[n]!$ is the rational function
$$
[n]! = \frac{\Gamma(\t z+n+1)}{\Gamma(\t z+1)} =
(1+z\t)(2+z\t)\dots(n+z\t) .
$$
The important cases for us will be
\begin{align*}
\pi_{-1}(z) &= \t + \lim_{\eps\to0} \sum_{n=1}^\infty
\frac{z^np_k(n)}{[n]!} , &
\pi_0(z) &= 1 + \lim_{\eps\to0} \sum_{n=1}^\infty \frac{z^np_0(n)}{[n]!} , &
\pi_1(z) &= \lim_{\eps\to0} \sum_{n=1}^\infty \frac{z^np_1(n)}{[n]!} .
\end{align*}
It follows from the recursion \eqref{recurse} for the coefficients of
$\LL$ that
\begin{equation} \label{recurse0}
\pi_{-1}(z) = q\pi_1(z)+\t \pi_0(z) .
\end{equation}
\begin{lemma} \label{epi}
$e(\pi_k(z)) = z\pi_k(z)$ and $\p_q\pi_k(z) = z\pi_{k+1}(z)$.
\end{lemma}
\begin{proof}
We have
\begin{align*}
(\LL_0-\t) e(\LL_0(z)) &= \sum_{n=-\infty}^\infty
\frac{nz^n}{[n]!} \LL_0^{n-1}(\LL_0-\t)e(\LL_0)
= \sum_{n=-\infty}^\infty \frac{nz^n}{[n]!} \LL_0^n \\
&= \sum_{n=-\infty}^\infty \frac{(\t z+n)z^n}{[n]!} \LL_0^n -
\sum_{n=-\infty}^\infty \frac{\t z^{n+1}}{[n]!} \LL_0^n =
z(\LL_0-\t)\LL_0(z) .
\end{align*}
The equation for $e(\pi_k(z))$ follows on taking the coefficient of
$\Lambda^k$. Since $\Lambda\p_qK_0=e(K_0)$, it follows from Theorem
\ref{LL0} that $\Lambda\p_q\LL_0(z)=e(\LL_0(z))$. Taking the coefficient of
$\Lambda^{k+1}$, we obtain the formula for $\p_q\pi_k(z)$.
\end{proof}

We may now prove the formulas \eqref{conjecture} and \eqref{conjecturebar}
relating the dispersionless limit of the equivariant Toda lattice to the
equivariant genus $0$ Gromov-Witten potential of $\CP^1$.

In the genus $0$ limit, the functions $q=\exp(u)$ and $v$ on the large
phase space become $\exp(\p^2\CF_0)$ and $\p\p_0\CF_0$, and $\vbar$ becomes
$v-\t u$. The proof of Theorem 4.2 of \cite{toda} extends to the
equivariant case, and shows that
\begin{align} \label{jet}
\p^nu &= s_n + O(|s|^2+|t|^2) , & \p^nv &= \delta_{1,n} + t_n +
O(|s|^2+|t|^2) .
\end{align}
Hence, we may identify the large phase space with the space of formal
jets in an affine space with coordinates $u$ and $v$.

The following lemma shows that the vector field $e$ lifts to the puncture
vector field on the large phase space.
\begin{lemma}
The puncture vector field
$$
e = \p - \sum_{n=0}^\infty \biggl( s_{n+1} \frac{\p}{\p s_n} + t_{n+1}
\frac{\p}{\p t_n} \biggr)
$$
on the large phase space acts on elements of $\Q[\t,q,v]$ by the derivation
$\p_v$.
\end{lemma}
\begin{proof}
Observe that the puncture vector field $e$ commutes with $\p$1; this
reflects the fact that $\CP^1$ is one-dimensional. The puncture (or string)
equation says that
\begin{equation} \label{string}
e(\CF_0)=s_0t_0+ \half\t t_0^2 .
\end{equation}
Applying the differential operators $\p_0\p$ and $\p^2$ to this
equation, we see that $e(v)=1$ and $e(u)=0$.
\end{proof}

\begin{theorem} \label{dispersionless}
The dispersionless limits of \eqref{conjecture} and
\eqref{conjecturebar} hold.
\end{theorem}
\begin{proof}
We will concentrate on the proof of \eqref{conjecture}. The proof of
the dispersionless limit of \eqref{conjecturebar} is the same, up to
conjugation.

Let $\p(z)$ be the generating function for vector fields
$$
\p(z) = \sum_{k=0}^\infty z^k \p_k .
$$
We must prove that
\begin{align*}
\p(z)v &= \sum_{n=1}^\infty \frac{z^{n-1}\p p_0(n)}{[n]!} , &
\p(z)u &= \sum_{n=1}^\infty \frac{z^{n-1}\p p_{-1}(n)}{[n]!} .
\end{align*}
In terms of the generating functions
\begin{align*}
x &= 1+z\p(z)\p\CF_0 - g_0(z) , &
y &= \t+z\p(z)\p_0\CF_0 - g_{-1}(z) ,
\end{align*}
we wish to prove that $\p x(z)=\p y(z)=0$. We will actually prove the
stronger result, that $x(z)=y(z)=0$: in other words, that
\begin{align*}
\sum_{k=0}^\infty z^k \p_k\p\CF_0 &= \sum_{n=1}^\infty
\frac{z^{n-1}p_0(n)}{[n]!} , &
\sum_{k=0}^\infty z^k \p_k\p_0\CF_0 &= \sum_{n=1}^\infty
\frac{z^{n-1}p_{-1}(n)}{[n]!} .
\end{align*}

A theorem of Dijkgraaf and Witten \cite{DW} establishes that the Toda
equation \eqref{Toda} holds in the dispersionless limit:
$$
\p_0^2\CF_0 = q + \t v .
$$
Combining the topological recursion relations for
equivariant Gromov-Witten invariants in genus $0$ with Lemma
\ref{epi}, we see that
\begin{align*}
\p x(z) &= z \bigl( x(z) \, \p(v-\t u) + y(z) \, \p u
\bigr) \\
\p y(z) &= z \bigl( x(z) \, \p q + y(z) \, \p v \bigr) .
\end{align*}
On the other hand, the string equation shows that $e(x(z)) =
zx(z)$ and $e(y(z)) = zy(z)$.

Now apply the following principle (Proposition 4.1 of \cite{toda}):
$$
\emph{A function $f$ on the large phase space such that $\p f$
and $e(f)$ lie in $\Q[\t]$ itself lies in $\Q[\t]$.}
$$
Arguing by induction, we see that the coefficients of $z^k$ in
$x(z)$ and $y(z)$ lie in $\Q[\t]$; in other words,
$x(z),y(z)\in\Q\[\t,z\]$. (In particular, we see that
$\p_k\p\CF_0$ and $\p_k\p_0\CF_0$ lie in $\Q[\t,q,v]$ for all
$k\ge0$.)

The proof is finished by observing that, by the divisor equation for
Gromov-Witten invariants, the limits $\lim_{q\to0}x(z)$ and
$\lim_{q\to0}y(z)$ are integrals over the degree $0$ moduli space
$\Mbar_{0,2}(\CP^1,0)$; however, this moduli space is empty, hence
$x(z)=y(z)=0$.
\end{proof}

As mentioned in the introduction, the analogue of Theorem
\ref{dispersionless} is now known to hold in all genera (Okounkov and
Pandharipande \cite{OP}).

\end{document}